\newcommand{\al}{\alpha}               
\newcommand{\be}{\beta}
        
\newcommand{\vphi}{\varphi}
\newcommand{\ga}{\gamma}               

\newcommand{\de}{\delta}

\newcommand{\sig}{\sigma}

\newcommand{\cal}{\mathcal}

\newcommand{\calf}{{\cal F}}

\newcommand{\calv}{{\cal V}}

\newcommand{\limpl}{\Longrightarrow}

\newcommand{\oo}{\infty}

\newcommand{\dst}{\displaystyle}

\def\R+oo{R_+\cup\{\oo\}}

\def\dtends  {\stackrel {\it d}{\longrightarrow}}

\def\(V)tends  {\stackrel {(\calv)}{\longrightarrow}}

\newcommand{\barr}{\begin{array}}       
\newcommand{\earr}{\end{array}}
\newcommand{\bcor}{\begin{corollary}}   
\newcommand{\ecor}{\end{corollary}}
\newcommand{\ben}{\begin{enumerate}}      
\newcommand{\een}{\end{enumerate}}
\newcommand{\beq}{\begin{equation}}     
\newcommand{\eeq}{\end{equation}}
\newcommand{\bit}{\begin{itemize}}      
\newcommand{\eit}{\end{itemize}}
\newcommand{\blemma}{\begin{lemma}}     
\newcommand{\elemma}{\end{lemma}}
\newcommand{\bproof}{\begin{proof}}     
\newcommand{\eproof}{\end{proof}}
\newcommand{\bprop}{\begin{proposition}} 
\newcommand{\eprop}{\end{proposition}}
\newcommand{\brem}{\begin{remark}}      
\newcommand{\erem}{\end{remark}}
\newcommand{\btab}{\begin{tabular}}     
\newcommand{\etab}{\end{tabular}}
\newcommand{\btheorem}{\begin{theorem}} 
\newcommand{\etheorem}{\end{theorem}}

\documentclass[reqno]{amsart}

\newtheorem{theorem}{\bf Theorem}
\newtheorem{corollary}{\bf Corollary}

\newtheorem{lemma}{\bf Lemma}
\newtheorem{proposition}{\bf Proposition}
\newtheorem{remark}{\bf Remark}

\begin{document}

\title
[Banach-Caristi contractive Pairs]
{COMMON  FIXED POINTS FOR \\
BANACH-CARISTI CONTRACTIVE PAIRS}

\author{Mihai Turinici}
\address{
"A. Myller" Mathematical Seminar;
"A. I. Cuza" University;
700506 Ia\c{s}i, Romania
}
\email{mturi@uaic.ro}


\subjclass[2010]{
54H25 (Primary), 06A06 (Secondary).
}

\keywords{
Metric space, Banach-Caristi contractive pair, common fixed point.
}

\begin{abstract}
Further extensions are given for the common fixed point statement in
Dien [J. Math. Anal. Appl., 187 (1994), 76-90],
involving Banach-Caristi contractive pairs.
\end{abstract}

\maketitle

\section{Introduction}
\setcounter{equation}{0}

Let $(X,d)$ be a complete metric space;
and $T\in \calf(X)$ be a selfmap of $X$.
[Here, for each couple of sets 
$\{A,B\}$, $\calf(A,B)$ stands for the class of all
functions from $A$ to $B$;
if $A=B$, one writes $\calf(A,A)$  as $\calf(A)$].
The following fixed point statement in
Caristi and Kirk  \cite{caristi-kirk-1975}
(referred to as: Caristi-Kirk theorem)
is our starting point.

\btheorem \label{t1}
Assume that
there exists a function $\al:X\to R_+:=[0,\oo[$ with
\ben
\item[] (a01)\ \ 
$d(x,Tx)\le \al(x)-\al(Tx)$,\ for each $x\in X$
\item[] (a02)\ \ 
$\al(.)$ is lsc on $X$\ ($\liminf_n \al(x_n)\ge \al(x)$,
whenever $x_n \dtends x$).
\een
Then, $T$ has at least one fixed point in $X$.
\etheorem

Note that, in terms of the [associated to $\al(.)$] order
\ben
\item[] (a03)\ \ 
($x,y\in X$):\ $x\le y$ iff $d(x,y)\le \al(x)-\al(y)$
\een
the contractive condition (a01) becomes
\ben
\item[] (a04)\ \ 
$x\le Tx$,\ for each $x\in X$\
(i.e.:\ $T$ is {\it progressive} on $X$).
\een
So, Theorem \ref{t1} is deductible from  
the arguments in 
Ekeland's variational principle \cite{ekeland-1979}.
Further aspects may be found in 
Brezis and Browder \cite{brezis-browder-1976};
see also
Turinici \cite{turinici-2009}.

Now, the Caristi-Kirk theorem 
found (especially via Ekeland's approach)
some basic applications to control and optimization, generalized
differential calculus, critical point theory  and normal solvability;
see the above references for details.
As a consequence,
many extensions of this result were proposed.
Here, we shall concentrate on the 1981 statement in this area
due to
Bhakta and Basu \cite{bhakta-basu-1981}.
Let $\{S,T\}$ be a couple of selfmaps in $\calf(X)$.
We say that $z\in X$ is a {\it common fixed point} of
$\{S,T\}$ if $Sz=Tz=z$.
Sufficient conditions guaranteeing such a property
are obtainable as below.
Call the selfmap $U$ in $\calf(X)$,
{\it orbitally continuous} (on $X$), if
\ben
\item[] (a05)\ \ 
$z=\lim_nU^{i(n)}x$ implies $Uz=\lim_n U^{i(n)+1}x$;
\een
here, $(i(n); n\ge 0)$ is a rank sequence with 
$i(n) \to \oo$ as $n\to \oo$.

\btheorem \label{t2}
Suppose that 
\ben
\item[] (a06)\ \ 
both $S$ and $T$ are orbitally continuous;
\een
and let the functions
$\al,\be:X\to R_+$ be such that
\ben
\item[] (a07)\ \ 
$d(Sx,Ty)\le \al(x)-\al(Sx)+\be(y)-\be(Ty)$,\
for all $x,y\in X$.
\een
Then,

{\bf i)} $S$ and $T$ have a unique common fixed point  $z\in X$,

{\bf ii)} $S^nx\dtends z$ and $T^nx\dtends z$ as $n\to \oo$, for each $x\in X$.
\etheorem

A partial extension of this result was given in the 1994 paper by
Dien \cite{dien-1994}:

\btheorem \label{t3}

Suppose that (a06) holds.
In addition, let
the number $q\in [0,1[$ and  
the function $\al:X\to R_+$ be such that
\ben
\item[] (a08)\ \ 
$d(Sx,Ty)\le qd(x,y)+\al(x)-\al(Sx)+\al(y)-\al(Ty)$,\
$\forall x,y\in X$.
\een
Then, conclusions of Theorem \ref{t2} are retainable.
\etheorem

[As a matter of fact, the original result is with
$\al=\al_1+...+\al_k$, where $\{\al_i; 1\le i\le k\}$
is a finite system in $\calf(X,R_+)$.
But it gives, practically, the same amount
of information as the result in question].

Note that, when $\al(.)$ is a constant function and
$S=T$, then Theorem \ref{t3} implies the 
Banach contraction principle \cite{banach-1922}.
In addition, (a01)  follows from (a08) when 
$S=I$ (=the {\it identity map} of $\calf(X)$)
and $x=y$;  for this reason, the couple $\{S,T\}$ above 
will be referred to as  {\it Banach-Caristi contractive}.
It is to be stressed that
Theorem \ref{t1} does not follow from Theorem \ref{t3};
because, the (essential for  Theorem \ref{t1}) condition (a02) 
is not obtainable from the conditions of Theorem \ref{t3}.
However, the underlying relationship between these results  holds whenever 
(a06) is accepted, in place of (a02).
(This clarifies an assertion made in
Ume and Yi \cite{ume-yi-2002};
we do not give details).
On the other hand, Dien's result cannot be deduced
from Caristi-Kirk's; because (a06) cannot be deduced from 
the conditions of Theorem \ref{t1}.
Finally, Theorem \ref{t2} cannot be viewed as a
particular case of Theorem \ref{t3}
(when $q=0$); because
the functions $\al(.)$ and $\be(.)$ may be distinct.

Having this precise, it is our aim in the following
to establish 
(cf.  Section 2) 
a common extension of these statements;
as well as (in  Section 3) 
a sum approach of it.
Some other aspects will be delineated elsewhere.

\section{Main result}
\setcounter{equation}{0}

Let $\vphi\in \calf(R_+)$ be a function;
call it {\it regressive} provided
$\vphi(0)=0$ and $\vphi(t)< t$, $\forall t\in R_+^0:=]0,\oo[$;
the class of all these will be 
denoted as $\calf(re)(R_+)$.
For example, any function $\vphi=qJ$ 
where $q\in [0,1[$, is regressive;
here, $J$ is the {\it identity function} of $\calf(R_+)$
($J(t)=t$, $t\in R_+$).

Now, fix some $\vphi\in \calf(re)(R_+)$.
Denote $\psi:=J-\vphi$; 
and call it, the {\it complement} of $\vphi$.
Clearly, $\psi\in \calf(R_+)$; precisely,
\beq \label{201}
\psi(0)=0;\ 0< \psi(t)\le t,\ \forall t\in R_+^0.
\eeq
For an easy reference, we list our basic hypotheses.
The former of these is
\ben
\item[] (b01)\ \ 
$\vphi$ is super-additive:\
$\vphi(t+s)\ge \vphi(t)+\vphi(s)$, for all $t,s\ge 0$;
\een
clearly, $\vphi$ must be increasing in such a case.
And the latter condition writes:
\ben
\item[] (b02)\ \ 
$\psi:=J-\vphi$ is coercive:\ 
$\psi(t)\to \oo$ as $t \to \oo$;
\een
referred to as:
$\vphi$ is {\it complementary coercive}.
Note that (by this very definition) 
\beq \label{202}
g(r):=\sup\{t\ge 0; \psi(t)\le r\}< \oo,\ \
\mbox{for each}\ r\in R_+;
\eeq
whence, $g(.)$ is an element of $\calf(R_+)$.
Moreover (from (\ref{201}) above)
\beq \label{203}
g(0)=0;\ g(r)\ge r,\ \forall r\in R_+.
\eeq

The following auxiliary fact will be useful for us.

\blemma \label{le1}
Let $\vphi\in \calf(re)(R_+)$ be super-additive and
complementary coercive. Further, 
let the sequence $(\theta_n; n\ge 0)$ in $R_+$ be such that
\ben
\item[] (b03)\ \ 
$\theta_{m+1}\le \vphi(\theta_m)+\de_m-\de_{m+1}$,\
for all $m\ge 0$;
\een
where $(\de_n; n\ge 0)$ is a sequence in $R_+$.
Then, the series $\sum_n\theta_n$ converges.
\elemma

\bproof
Let $(\sig_i:=\theta_0+...+\theta_i; i\ge 0)$ be the
partial sum sequence attached to $(\theta_n; n\ge 0)$.
For each $j\ge 0$, we have 
(summing in (b03) from $m=0$ to $m=j$)
$$
\theta_1+...+\theta_{j+1}\le
\vphi(\theta_0)+...+\vphi(\theta_j)+\de_0-\de_{j+1}.
$$
This, along with the super-additivity of $\vphi$, gives
$\sig_j\le \vphi(\sig_j)+\theta_0+\de_0$.
But then, (\ref{202}) yields
[$\sig_n \le g(\theta_0+\de_0)< \oo$, $\forall n$];
wherefrom, all is clear.
\eproof

We now state the promised result.
Let $(X,d)$ be a complete metric space;
and $\{S,T\}$ be a pair in $\calf(X)$.

\btheorem \label{t4}
Suppose that (a06) holds.
In addition, let the function
$\vphi\in \calf(re)(R_+)$ as in (b01)+(b02), and the map
$\ga:X\times X \to R_+$, be taken so as
\ben
\item[] (b04)\ \ 
$d(Sx,Ty)\le \vphi(d(x,y))+\ga(x,y)-\ga(Sx,Ty)$,\
for all $x,y\in X$.
\een
Then, conclusions of Theorem \ref{t2} are retainable.
\etheorem

\bproof
Given $x_0,y_0\in X$,  
put $(x_n=S^nx_0; n\ge 0)$, $(y_n=T^ny_0; n\ge 0)$.
From (b04), one has (by a finite induction argument),
the iterative type relations
$$
d(x_{m+1},y_{m+1})\le \vphi(d(x_m,y_m))+
\ga(x_m,y_m)-\ga(x_{m+1},y_{m+1}),\
\mbox{for all}\ m\ge 0.
$$
Combining with Lemma \ref{le1} (and the adopted notations),
one derives that
the series $\sum_nd(S^nx_0,T^ny_0)$ converges.

Further, let us develop the same reasoning by starting
from the points $u_0=Sx_0$ and $y_0$;
one derives that
the series $\sum_n d(S^nu_0,T^ny_0)$ converges;
or, equivalently:
the series $\sum_nd(S^{n+1}x_0,T^ny_0)$ converges.
This, along with
$$
d(S^nx_0,S^{n+1}x_0)\le d(S^nx_0,T^nx_0)+
d(S^{n+1}x_0,T^ny_0),\ \ \forall n\ge 0
$$
tells us that 
the series $\sum_n d(S^nx_0,S^{n+1}x_0)$ converges; 
wherefrom, 
$(S^nx_0; n\ge 0)$ is a $d$-Cauchy sequence.
In a similar way
(starting from the points $x_0$ and $v_0=Ty_0$) one proves that
the series $\sum_n d(T^ny_0,T^{n+1}y_0)$ converges; 
whence, 
$(T^ny_0; n\ge 0)$ is a $d$-Cauchy sequence.
As $(X,d)$ is complete, we have that 
$S^nx_0\dtends z$ and $T^ny_0\dtends w$, 
for some $z,w \in X$.
Combining with the orbital continuity of both $S$ and $T$,
gives 
$S(S^nx_0)\dtends Sz$, $T(T^ny_0)\dtends Tw$.
But, 
$S(S^nx_0)=S^{n+1}x_0\dtends z$, $T(T^ny_0)=T^{n+1}y_0\dtends w$;
and this yields $z=Sz$, $w=Tw$.
Finally, from (b04) again, we have $d(z,w)\le \vphi(d(z,w))$;
so that, $z=w$.
Hence, $z$ is a common fixed point of $\{S,T\}$.
Its uniqueness is obtainable 
by the argument we just developed for $(z,w)$;
and, from this, we are done.
\eproof

In particular, when $\vphi\in \calf(re)(R_+)$ is taken as
\ben
\item[] (b05)\ \ 
$\vphi(t)=qt$,\ $t\ge 0$, for some $q\in [0,1[$,
\een
conditions (b01)+(b02) hold; and then, under the choice 
\ben
\item[] (b06)\ \ 
$\ga(x,y)=\al(x)+\al(y)$, $x,y\in X$\ 
(where $\al\in \calf(X,R_+)$)
\een
the corresponding version of Theorem \ref{t4}
 is just Theorem \ref{t3}.
Note that, under the same framework, 
a more general choice for $\ga$ is
\ben
\item[] (b07)\ \ 
$\ga(x,y)=\al(x)+\be(y)$, $x,y\in X$\ 
(where $\al,\be\in \calf(X,R_+)$).
\een
This version of Theorem \ref{t4} is (under $\vphi=0$)
just Theorem \ref{t2} above. 
Further aspects may be found in 
Alimohammady et al
\cite{alimohammady-balooee-radojevic-rakocevic-roohi-2011};
see also 
Kadelburg et al \cite{kadelburg-radenovic-simic-2011}.

\section{Further extensions}
\setcounter{equation}{0}

A simple inspection of the argument we just developed 
shows that it depends (via Lemma \ref{le1}) on the
super-additivity of the function $\vphi\in  \calf(re)(R_+)$;
so, we may ask whether this cannot be removed.
An appropriate answer is available,
if we arrange for the sums given by the argument
of Theorem \ref{t4} being obtainable in a direct way
from the contractive conditions.

Let $(X,d)$ be a complete metric space;
and $\{S,T\}$ be a pair of selfmaps in $\calf(X)$.

\btheorem \label{t5}
Suppose that (a06) holds.
In addition, let the function
$\vphi\in \calf(re)(R_+)$ as in (b02) and the map
$\ga:X\times X\to R_+$ be taken so as
\ben
\item[] (c01)\ \ 
$\dst \sum_{j=1}^n d(S^jx,T^jy)\le
\vphi(\sum_{j=0}^{n-1} d(S^jx,T^jy))+
\ga(x,y)-\ga(S^nx,T^ny)$,
\een
for all $x,y\in X$ and all  $n\ge 1$.
Then, conclusions of Theorem \ref{t3} are retainable.
\etheorem

\bproof
Given $x_0,y_0\in X$, 
put $(x_n=S^nx_0; n\ge 0)$, $(y_n=T^ny_0; n\ge 0)$.
Further, denote $(\theta_i=d(x_i,y_i); i\ge 0)$.
By (c01) one has, for each $n\ge 1$,
$$
\theta_1+...+\theta_n\le
\vphi(\theta_0+...+\theta_{n-1})+ \ga(x_0,y_0)-\ga(x_n,y_n);
$$
wherefrom (after some transformations)
$$
\theta_0+...+\theta_{n-1}\le
\vphi(\theta_0+...+\theta_{n-1})+ \theta_0+\ga(x_0,y_0),\ 
\forall n\ge 1.
$$
This, from (b02) (and the notations in Section 2), gives
$$
\theta_0+...+\theta_{n-1}\le g[\theta_0+\ga(x_0,y_0)]< \oo,\
\mbox{for all}\ n\ge 1;
$$
so that (by the adopted notations)
the series $\sum_n d(S^nx_0,T^ny_0)$ converges.
Further, let us develop the same reasoning by starting
from the points $u_0=Sx_0$ and $y_0$;
one derives that
the series $\sum_nd(S^nu_0,T^ny_0)$ converges;
or, equivalently:
the series $\sum_nd(S^{n+1}x_0,T^ny_0)$ converges.
This, along with 
$$
d(S^nx_0,S^{n+1}x_0)\le d(S^nx_0,T^nx_0)+
d(S^{n+1}x_0,T^ny_0),\ \ \forall n\ge 0
$$
tells us that the series $\sum_nd(S^nx_0,S^{n+1}x_0)$
converges; wherefrom $(S^nx_0; n\ge 0)$ is a
$d$-Cauchy sequence.
In a similar way
(starting from the points $x_0$ and $v_0=Ty_0$), one proves that
the series $\sum_n d(T^ny_0,T^{n+1}y_0)$ converges; 
whence, 
$(T^ny_0; n\ge 0)$ is a $d$-Cauchy sequence.
As $(X,d)$ is complete, 
$S^nx_0\dtends z$ and $T^ny_0\dtends w$, 
for some $z,w \in X$.
The remaining part of the argument 
runs as in Theorem \ref{t4}, because (c01) $\limpl$ (b04); 
and, from this, all is clear.
\eproof

Now, concrete examples of 
complementary coercive functions 
$\vphi\in \calf(re)(R_+)$ are 
obtainable by starting from the choice
\ben
\item[] (c02)\ \ 
$\vphi(t)=t\chi(t)$,\ $t\ge 0$,
\een
where $\chi\in \calf(R_+)$ fulfills
the regularity conditions
\ben
\item[] (c03)\ \ 
$\chi$ is increasing on $R_+^0$ and $\chi(t)< 1, \forall t\in R_+^0$.
\een
The standard case is [$\chi(t)=q$, $t\ge 0$], 
for some $q\in [0,1[$; when Theorem \ref{t5} includes, in a
direct way, Theorem \ref{t3}; due, as above said, to
Dien \cite{dien-1994}.
A technical extension of this one may be constructed according to
\ben
\item[] (c04)\ \ 
$\chi(t)=r_{n+1}$,\ when $t\in [t_n,t_{n+1}[$, for each $n\ge 0$;
\een
where the sequence $(r_n; n\ge 1)$  in $]0,1[$
and the strictly ascending sequence 
$(t_n; n\ge 0)$ in $R_+$ with $t_0=0$
and $t_n\to \oo$ are to be determined.
To this end, we have 
$$
t-\vphi(t)=t(1-r_{n+1}),\ t\in [t_n,t_{n+1}[,\ n\ge 0.
$$
Assume that $(t_n; n\ge 1)$ 
is a strictly ascending sequence in $]1,\oo[$ with
\ben
\item[] (c05)\ \ 
$t_n/\sqrt{t_{n+1}} \to \oo$\ (hence $t_n\to \oo$).
\een
Then, choose the sequence $(r_n; n\ge 1)$ in $]0,1[$
according to
\ben
\item[] (c06)\ \ 
$1-r_n=1/\sqrt{t_n}$, for each $n\ge 1$.
\een
Note that, as a consequence of this, $(r_n; n\ge 1)$
is strictly ascending in $]0,1[$ (hence, (c03) holds)
and $r_n\to 1$ as $n\to \oo$.
Replacing in a preceding formula, yields
$$
t-\vphi(t)=t/\sqrt{t_{n+1}},\ \
\mbox{when}\ t\in [t_n,t_{n+1}[,\ n\ge 0.
$$
This gives an evaluation like
$$
t-\vphi(t)\ge t_n/\sqrt{t_{n+1}},\
\mbox{for}\ t\in [t_n,t_{n+1}[, n\ge 0;
$$
wherefrom (by (c05)), 
$\psi:=J-\vphi$ is coercive.
On the other hand, some useful 
super-additivity tests 
for the functions 
$\vphi\in \calf(re)(R_+)$ like before
are obtainable from the methods 
developed by 
Bruckner \cite{bruckner-1962}.
Some other aspects 
may be found in 
Liu, Xu and Cho \cite{liu-xu-cho-1998};
see also
Fisher \cite{fisher-1981}.


\end{document}